\documentclass[12pt]{article}
\usepackage{graphicx}
\usepackage{amsmath,amsthm,amssymb,enumerate}
\usepackage{euscript,mathrsfs}
\usepackage[left=2.5cm,right=2.5cm,top=3cm,bottom=3cm]{geometry}
\usepackage{color}
\catcode`\@=11 \@addtoreset{equation}{section}

\catcode`\@=12

\allowdisplaybreaks

\newtheorem{Theorem}{Theorem}[section]
\newtheorem{Proposition}[Theorem]{Proposition}
\newtheorem{Lemma}[Theorem]{Lemma}
\newtheorem{Corollary}[Theorem]{Corollary}

\theoremstyle{definition}
\newtheorem{Definition}[Theorem]{Definition}

\newtheorem{Remark}[Theorem]{Remark}

\newcommand{\bTheorem}[1]{
\begin{Theorem} \label{T#1} }
\newcommand{\eT}{\end{Theorem}}

\newcommand{\bProposition}[1]{
\begin{Proposition} \label{P#1}}
\newcommand{\eP}{\end{Proposition}}

\newcommand{\bLemma}[1]{
\begin{Lemma} \label{L#1} }
\newcommand{\eL}{\end{Lemma}}

\newcommand{\bCorollary}[1]{
\begin{Corollary} \label{C#1} }
\newcommand{\eC}{\end{Corollary}}

\newcommand{\bRemark}[1]{
\begin{Remark} \label{R#1} }
\newcommand{\eR}{\end{Remark}}

\newcommand{\bDefinition}[1]{
\begin{Definition} \label{D#1} }
\newcommand{\eD}{\end{Definition}}

\newcommand{\bu}{\mathbf u}

\newcommand{\bfphi}{\boldsymbol{\varphi}}

\newcommand{\bFormula}[1]{
\begin{equation} \label{#1}}
\newcommand{\eF}{\end{equation}}

\newcommand{\Ov}[1]{\overline{#1}}

\newcommand{\DC}{C^\infty_c}

\newcommand{\vr}{\varrho}
\newcommand{\vre}{\vr_\ep}

\newcommand{\vte}{\vt_\ep}
\newcommand{\vue}{\vu_\ep}
\newcommand{\tvr}{\tilde \vr}
\newcommand{\tvu}{{\tilde \vu}}
\newcommand{\tvt}{\tilde \vt}

\newcommand{\vt}{\vartheta}
\newcommand{\vu}{\vc{u}}
\newcommand{\vm}{\vc{m}}
\newcommand{\Ee}{E_{\ep}}

\newcommand{\vq}{\vc{q}}
\newcommand{\vn}{\vc{n}}
\newcommand{\vme}{\vm_\ep}
\newcommand{\vc}[1]{{\bf #1}}

\newcommand{\Div}{{\rm div}_x}
\newcommand{\Grad}{\nabla_x}
\newcommand{\Dt}{\frac{\rm d}{{\rm d}t}}

\newcommand{\dx}{\,{\rm d} {x}}

\newcommand{\dt}{\,{\rm d} t }

\newcommand{\intO}[1]{\int_{\Omega} #1 \ \dx}

\newcommand{\ep}{\varepsilon}

\newcommand{\R}{\mathbb{R}}
\newcommand{\I}{\mathbb{I}}
\renewcommand{\S}{\mathbb{S}}

\definecolor{Cgrey}{rgb}{0.85,0.85,0.85}
\definecolor{Cblue}{rgb}{0.50,0.85,0.85}
\definecolor{Cred}{rgb}{1,0,0}
\definecolor{fancy}{rgb}{0.10,0.85,0.10}

\newcommand\Cbox[2]{%
    \newbox\contentbox%
    \newbox\bkgdbox%
    \setbox\contentbox\hbox to \hsize{%
        \vtop{
            \kern\columnsep
            \hbox to \hsize{%
                \kern\columnsep%
                \advance\hsize by -2\columnsep%
                \setlength{\textwidth}{\hsize}%
                \vbox{
                    \parskip=\baselineskip
                    \parindent=0bp
                    #2
                }%
                \kern\columnsep%
            }%
            \kern\columnsep%
        }%
    }%
    \setbox\bkgdbox\vbox{
        \color{#1}
        \hrule width  \wd\contentbox %
               height \ht\contentbox %
               depth  \dp\contentbox
        \color{black}
    }%
    \wd\bkgdbox=0bp%
    \vbox{\hbox to \hsize{\box\bkgdbox\box\contentbox}}%
    \vskip\baselineskip%
}

\date{}

\begin{document}


\title{Existence of measure-valued solutions to a complete Euler system for a perfect gas}

\author{Jan B\v rezina}

\date{\today}

\maketitle

\bigskip

\centerline{Tokyo Institute of Technology}

\centerline{ 2-12-1 Ookayama, Meguro-ku, Tokyo, 152-8550, Japan}

\bigskip

\begin{abstract}
The concept of {\it renormalized dissipative measures-valued} (rDMV) {\it solutions} to a {\it complete Euler system} for a {\it perfect gas} was introduced in  \cite{BreFei17} and further discussed in \cite{BreFei17A}. Moreover it was shown there that rDMV solutions satisfy the weak (measure--valued)--strong uniqueness principle that makes them a useful tool.  In this paper we prove the existence of rDMV solutions.
Namely,  we formulate the complete Euler system in conservative variables usual for numerical analysis and recall the concept of rDMV solutions based on the total energy balance and renormalization of entropy inequality for the physical entropy presented in \cite{BreFei17}. We then give two different ways how to generate rDMV solutions. First  via vanishing viscosity limit using Navier-Stokes equations coupled with entropy transport and second via the vanishing dissipation limit of the two-velocity model proposed by H. Brenner. Finally, we recall  the weak--strong uniqueness principle for rDMV solutions proved in \cite{BreFei17} and \cite{BreFei17A}.

\end{abstract}



\section{Introduction}

In this paper we show how to generate {\it renormalized dissipative measures-valued} (rDMV) {\it solutions} to a {\it complete Euler system} for a {\it perfect gas} by two different ways.  Namely, we  are interested in the following problem.

Let $\Omega \subset R^N$ be a physical domain occupied by the fluid. For simplicity we assume that $N=3$ and $\Omega$ is a bounded domain with a smooth boundary $\partial \Omega$. 
The time evolution of the mass density $\vr = \vr(t,x)$, the velocity field $\vu = \vu(t,x)$, and the (absolute) temperature $\vt = \vt(t,x)$ is governed by the following system of partial differential equations expressing the basic physical principles:

\begin{itemize}
\item {\bf Conservation of mass}
\begin{equation} \label{e1}
\partial_t \vr + \Div (\vr \vu) = 0;
\end{equation}

\item {\bf Conservation of linear momentum}
\begin{equation} \label{e2}
\partial_t (\vr \vu) + \Div (\vr \vu \otimes \vu) + \Grad (\vr \vt) = 0;
\end{equation}

\item {\bf Conservation of total energy}

\begin{equation} \label{e3}
\partial_t \left( \frac{1}{2} \vr |\vu|^2 + c_v\vr\vt \right) +
\Div \left[ \left( \frac{1}{2} \vr |\vu|^2 + c_v\vr\vt + \vr \vt \right) \vu \right] = 0,
\end{equation}
where $ c_v > 0 $ is the specific heat at constant volume.
\end{itemize}

Moreover, in accordance with  the {\bf Second law of thermodynamics}, the entropy $$s= s(\vr,\vt) = \log (\vt^{c_v}) - \log (\vr)$$ should satisfy the equation
\begin{equation} \label{e4}
\partial_t \left(\vr \log\left(\frac{\vt^{c_v}}{\vr}\right)\right) + \Div \left(\vr \log\left(\frac{\vt^{c_v}}{\vr}\right)\vu\right) = 0.
\end{equation}
It is easy to check that \eqref{e4} follows from (\ref{e1}--\ref{e3}) as long as all quantities in (\ref{e1}-\ref{e3}) are continuously differentiable.

As is well known, smooth solutions of (\ref{e1}-\ref{e3}) exist only for a finite lap of time after which singularities develop for a fairly generic class of initial data.
Therefore global-in-time solutions may exist only in a weak sense, where the derivatives in (\ref{e1}-\ref{e3}) are understood in the sense of distributions. In that case \eqref{e4} is no longer automatically satisfied and it may be added as an admissibility condition.

In view of recent results based on the theory of convex integration, see \cite{FeKlKrMa}, weak solutions of (\ref{e1}--\ref{e3}), even if supplemented by
(\ref{e4}), are not uniquely determined by the initial data
as long as $N > 1$. As a matter of fact, for any \emph{piecewise constant} initial density $\vr_0$ and temperature $\vt_0$, there exists
$\vu_0 \in L^\infty(\Omega; R^N)$, $N=2,3$ such that the problem (\ref{e1}--\ref{e4}) with (\ref{e7}) admits infinitely many weak (distributional) solutions on a given time interval $(0,T)$. This kind of result indicates that we should look for a different approach to   concept of solutions to the  Euler system.

In his pioneering work \cite{DiP2}, DiPerna proposed a new concept of solution, known as \emph{measure--valued solution}, to nonlinear systems of partial differential equations
admitting uncontrollable oscillations. In particular with focus on the compressible Euler system and other related models of \emph{inviscid} fluids.
Although \emph{existence} of a measure-valued solution to a given problem is usually an almost straightforward consequence of {\it a priori} bounds, its \emph{uniqueness} in terms of the initial data can be seen as the weakest point of this approach. On the other hand, Brenier et al. \cite{BrDeSz} proposed a new approach seeing the measure-valued solutions as possibly the largest class in which the family of  smooth (classical) solutions
is stable. In particular, they show the so-called weak (measure-valued)--strong uniqueness principle
for the incompressible Euler system. Specifically,
a classical and a measure--valued solution emanating from the same initial data coincide as long as the former exists. These results have been extended
to the isentropic Euler and Navier--Stokes systems by Gwiazda et al. \cite{GSWW} and \cite{FGSWW1}. 
Following the philosophy of Brenier et al. \cite{BrDeSz}, we focus on the concept of measure-valued solutions in the widest possible sense. Accordingly, using the fundamental laws of
thermodynamics, we extract the minimal piece of information to be retained to preserve the weak--strong uniqueness principle. To do so,  we follow the approach advocated in \cite{FENO6}, where equations (\ref{e1}), (\ref{e2}) and  (\ref{e4})  are
supplemented with the total energy inequality
\begin{equation} \label{e5}
\Dt \intO{ \left[ \frac{1}{2} \vr |\vu|^2 + \vr e( \vr, \vt) \right] } \leq 0.
\end{equation}

Moreover, since integrability of the convective term in the entropy equality (\ref{e4}) is problematic,  motivated by the work of Chen and Frid \cite{CheFr2} we consider a ``regularized'' version of  (\ref{e4}) relaxed to {\it inequality}, in particular, we consider
\begin{equation} \label{e6}
\partial_t \left(\vr \chi\left(\log\left(\frac{\vt^{c_v}}{\vr}\right)\right)\right) + \Div \left(\vr \chi\left( \log\left(\frac{\vt^{c_v}}{\vr}\right)\right)\vu\right) \geq 0,
\end{equation}
for any  increasing concave function $\chi$ satisfying $\chi(s) \leq \chi_\infty$ for all $s\in R$.

The problem is closed by prescribing the initial data and slip boundary conditions
\begin{equation} \label{e7}
\vr(0, \cdot) = \vr_0, \ \vt(0, \cdot) = \vt_0, \ \vu(0, \cdot) = \vu_0,
\end{equation}
\begin{equation} \label{eBC} 
\vu \cdot \vc{n}|_{\partial \Omega} = 0.
\end{equation}

The measure--valued solutions are natural candidates for describing the zero dissipation limits of more complex systems of Navier--Stokes type. At this point it is rather disappointing that we do not know how to construct rDMV solutions via some vanishing dissipation limit of Navier--Stokes--Fourier system even though we have such result for other types of fluids (see \cite {BreFei17A}). The issue is, as expected, to get good enough estimates to pass in the limit.

To overcome this difficulty we consider instead two other models to generate rDMV solutions.  Namely, Navier-Stokes with entropy transport and a two-velocity model proposed by H.Brenner \cite{BREN2}, \cite{BREN}, \cite{BREN1}. In particular it is interesting to see that Brenner's model behaves actually better in the vanishing dissipation limit and does not suffer the drawbacks of Navier-Stokes-Fourier system.

The paper is organized as follows. In Section \ref{MV}, we reformulate the problem in {\it conservative variables} and define rDMV  solutions with respect to the new formulation as in \cite{BreFei17} and \cite{BreFei17A}.  In Section \ref{ES}, we prove the existence of rDMV solutions via vanishing viscosity limit using  Navier-Stokes equations with variable entropy. Actually the solutions will satisfy \eqref{e6} as an equality even for any  continuous function $\chi$. 
In Section \ref{EB}, we show that rDMV solutions can be also generated via vanishing dissipation limit using Brenner's two velocity model.
Finally, in Section \ref{WS}, we state the weak (measure-valued)--strong uniqueness principle and the ideas behind its proof published in \cite{BreFei17} and \cite{BreFei17A}.

\section{Measure-valued solutions}
\label{MV}

\subsection{Conservative variables}

To introduce the concept of measure--valued solution for complete Euler system, it is more convenient to formulate the problem in the \emph{conservative variables}:
\[
\vr, \ \vc{m}= \vr \vu, \ E = \frac{1}{2} \frac{|\vc{m}|^2}{\vr} + c_v\vr \vt, \ c_v >0.
\]
The reason for changing the phase space is the fact that the temperature $\vt$ as well as the velocity
$\vu$ may not be correctly defined on the (hypothetical) vacuum set.
As the measure--valued solutions are typically (see DiPerna \cite{DiP2}) generated as weak limits of suitable approximation schemes, the presence of vacuum zones cannot be
{\it a priori} excluded. Moreover, such formulation is common in numerical analysis, where we aim to utilize our results.

The system (\ref{e1}), (\ref{e2}), (\ref{e5}), (\ref{e6}) with (\ref{eBC}) rewrites as
\begin{equation}
\label{p1}
\begin{split}
\partial_t \vr + \Div \vc{m} &= 0,\\
\partial_t \vc{m} + \Div \left( \frac{\vc{m} \otimes \vc{m}}{\vr} \right) + \frac1{c_v} \Grad \left( E - \frac{1}{2} \frac{|\vc{m}|^2}{\vr} \right)  &= 0,\\
\partial_t \intO{  E }  &\leq  0,
\end{split}
\end{equation}
together with the associated ``renormalized'' entropy inequality
\begin{equation} \label{p2}
\partial_t \left( \vr \chi\left(c_v \log \left( \ \frac{E - \frac{1}{2} \frac{|\vc{m}|^2}{\vr} }{c_v \vr^\gamma} \right)\right) \right)
+ \Div \left[ \chi\left(c_v \log \left(  \frac{E - \frac{1}{2} \frac{|\vc{m}|^2}{\vr} }{c_v \vr^\gamma} \right) \vc{m} \right)\right]  \geq 0,
\end{equation}
for any  increasing concave function $\chi$ satisfying $\chi(s) \leq \chi_\infty$ for all $s\in R$ and the slip boundary condition
\begin{equation} \label{pBC}
\vm \cdot \vc{n}|_{\partial \Omega} = 0.
\end{equation}

Although the thermodynamic functions are well defined for regular values $\vr > 0$, $\vt > 0$ of the standard variables, where the latter condition corresponds in the conservative setting
to $E - \frac{1}{2} \frac{|\vc{m}|^2}{\vr} > 0$, we need them to be defined even for the limit values $\vr = 0$, $\vt = 0$. To that end, we first define
\[
\frac{1}{2} \frac{|\vc{m}|^2}{\vr} = \left\{ \begin{array}{l} \frac{1}{2} \frac{|\vc{m}|^2 }{\vr} \ \mbox{for}\ \vr > 0,\\ \\
0 \ \mbox{if} \ \vc{m} = 0, \\ \\ \infty \ \mbox{otherwise.}
\end{array} \right.
\]
Like that we get  a lower semi--continuous convex function defined on the set $\{ \vr \geq 0, \ \vc{m} \in R^3 \}$. Second, we introduce the {\it renormalized total entropy}
$$
\mathcal{S}_\chi(\vr, \vc{m}, E)
=
\left\{ \begin{array}{l}
\vr \chi \left( c_v \log \left( \frac{E - \frac{1}{2} \frac{|\vc{m}|^2}{\vr} }{c_v \vr^\gamma} \right)\right) \ \mbox{if}\ \vr > 0,
\ E >  \frac{1}{2} \frac{|\vc{m}|^2}{\vr},\\ \\
0 \ \mbox{if}\ \vr = 0, \ \vc{m} = 0, \ E \geq 0,\\ \\
- \infty \ \mbox{otherwise}.
\end{array} \right.
$$
The total entropy $\mathcal{S}_\chi$ defined this way is a concave upper semi--continuous function defined on the set $\{ \vr \geq 0,\ \vc{m} \in R^3, \ E \geq 0 \}$ for every non-decreasing concave function $\chi$ satisfying $\chi(s) \leq \chi_\infty$ for all $s\in R$ (see  \cite{BreFei17B}).

\subsection{Renormalized dissipative measure--valued  solutions}

The following definition of the  renormalized dissipative measure--valued (rDMV) solutions was introduced in \cite{BreFei17}.

The initial state of the system is given through a parameterized family of probability measures $\{ U_{0,x} \}_{x \in \Omega}$ defined on the
phase space
\[
\mathcal{Q} \equiv
\left\{ (\vr, \vc{m}, E) \ \Big| \ \vr \geq 0, \ \vc{m} \in R^3, \ E \geq 0 \right\},
\]
andit is assumed
that the mapping $x \mapsto U_{0,x}$ belongs to $L^\infty_{{\rm weak-(*)}}(\Omega; \mathcal{P}(\mathcal{Q}))$.

Similarly, an rDMV solution is represented by a family of probability measures
\[
\{ U_{t,x} \}_{(t,x) \in (0,T) \times \Omega},\ U \in L^\infty_{{\rm weak}-(*)} ((0,T) \times \Omega; \mathcal{P}(\mathcal{Q}) ),
\]
and the non-linearities in (\ref{p1}), (\ref{p2}) are replaced by their expected values whereas the derivatives are understood in the sense of distributions. 

Hereafter $\left< U_{t,x}, g(\vr, \vc{m}, E) \right>$ denotes the expected value of a (Borel) function $g$ defined on $\mathcal{Q}$.

\begin{Definition} \label{D1}
A parameterized family of probability measures $U \in L^\infty_{{\rm weak}-(*)} ((0,T) \times \Omega; \mathcal{P}(\mathcal{Q}) )$ is
called a \emph{renormalized dissipative measure--valued} (rDMV) {\it solution} to the Euler system (\ref{p1}--\ref{pBC}) with the initial data
$U_0 \in L^\infty_{{\rm weak}-(*)} (\Omega; \mathcal{P}(\mathcal{Q}) )$ if the following holds:
\begin{itemize}
\item
\begin{equation} \label{p4}
\int_0^T \intO{ \left[ \left< U_{t,x}; \vr \right> \partial_t \varphi + \left< U_{t,x}; \vc{m} \right> \cdot
\Grad \varphi \right] } \dt = - \intO{ \left< U_{0,x} ; \vr \right> \varphi(0, \cdot)}
\end{equation}
for any $\varphi \in \DC([0,T) \times \Ov{\Omega})$;
\item
\begin{equation} \label{p5}
\begin{split}
\int_0^T & \intO{ \left[ \left< U_{t,x}; \vc{m} \right> \cdot \partial_t \bfphi  + \left< U_{t,x}; \frac{ \vc{m} \otimes \vc{m} }{\vr}
\right> : \Grad \bfphi + \frac1{c_v} \left< U_{t,x}; E - \frac{1}{2} \frac{|\vc{m}|^2 }{\vr} \right> \Div \bfphi \right] }\dt \\
&= - \intO{ \left< U_{0,x}; \vc{m} \right> \cdot \bfphi(0, \cdot) } +
\int_0^T \int_{\Ov{\Omega}} \Grad \bfphi : {\rm d} \mu_C
\end{split}
\end{equation}
for any $\bfphi \in \DC([0, T) \times \Ov{\Omega}; R^3)$, $ \bfphi\cdot \vn|_{\partial \Omega} = 0$  and  $\mu_C$ is a (vectorial) signed measure on $[0,T]\times \Ov{\Omega}$;
\item
\begin{equation} \label{p6}
\intO{ \left< U_{\tau,x}; E \right> } \leq \intO{ \left< U_{0,x}; E \right> } \ \mbox{for a.a.}\ \tau \in (0,T);
\end{equation}
\item
\begin{equation} \label{p7}
\begin{split}
\int_0^T &\intO{ \left[ \left< U_{t,x} ; \mathcal{S}_\chi (\vr, \vc{m}, E) \right> \partial_t \varphi +
\left< U_{t,x}; \mathcal{S}_\chi (\vr, \vc{m}, E) \frac{\vc{m}}{\vr} \right> \cdot \Grad \varphi  \right] } \dt \\
&\leq - \intO{ \left< U_{0,x}; \mathcal{S}_\chi(\vr, \vc{m}, E) \right> \varphi(0, \cdot) }
\end{split}
\end{equation}
for any $\varphi \in \DC([0,T) \times \Ov{\Omega})$, $\varphi \geq 0$, and any increasing concave function $\chi$ satisfying $\chi(s) \leq \chi_\infty$ for all $s\in R$;
\item
\begin{equation} \label{p8}
\int_0^\tau \int_\Omega d \left| \mu_C \right|  \leq c(c_v) \int_0^\tau \intO{ \left[ \left< U_{0,x}; E \right> - \left< U_{t,x}; E \right> \right]}
\dt \ \mbox{for any}\ 0 \leq \tau < T.
\end{equation}
\end{itemize}
\end{Definition}

Any ``standard'' weak solution $(\vr, \vc{m}, E)$ to (\ref{p1}--\ref{pBC}) may be identified with a measure--valued solution $U$ via
\[
U_{t,x} = \delta_{\vr(t,x), \vc{m}(t,x), E(t,x) } \ \mbox{for a.a.}\ (t,x) \in (0,T) \times \Omega,
\]
where $\delta_Z$ denotes the Dirac measure supported by $Z$.

The Definition \ref{D1} was motivated by previous works of others as well as by a result from \cite{BreFei17}. In particular, we showed in \cite{BreFei17} that any cluster point of a family  of "standard" admissible weak solutions to (\ref{e1}), (\ref{e2}), (\ref{e4}), (\ref{e6}) with uniformly bounded initial data is
an rDMV solution in the sense of Definition \ref{D1}. We point out that our class of measure-valued solutions includes all admissible weak solutions to the Euler system. Moreover, we showed that the rDMV solutions enjoy certain minimum principle out of the vacuum set, that is
\[
U_{0,x} \left\{ s(\vr, \vc{m},E) \geq s_0 \right\} = 1 \ \mbox{implies}\
U_{t,x} \left\{ s(\vr, \vc{m},E) \geq s_0\ | \ \vr >0 \right\} = 1 \ \mbox{for a.a.} \ (t,x).
\]

To conclude, we remark that the family of rDMV solutions for a given initial data is closed with respect to convex combinations. In particular, in view of the results obtained in \cite{FeKlKrMa}, there is a vast class of initial data for which the Euler system admits infinitely many nontrivial rDMV solutions. Here nontrivial means that  they do not consist of a single Dirac mass.

\section{Existence of rDMV solutions via Navier-Stokes}
\label{ES}

In this section we show the existence of rDMV solutions via the vanishing viscosity limit for a Navier-Stokes system with variable entropy
\begin{equation} \label{NSve}
\begin{split} 
 \partial_t \vr + \Div (\vr \vu)&= 0, 
\\
\partial_t (\vr \vu) + \Div \left(\vr \vu\otimes \vu\right) + \Grad p  &=  \Div \mathbb{S}(\nabla \vu),
\\
 \partial_t \left(\vr s\right) + \Div \left(\vr s\vu\right)  &= 0.
\end{split}
\end{equation}
Here, $p$ denotes the pressure; $\mathbb{S}$ is the viscous stress tensor 
$$
\mathbb{S}(\Grad \vu) = \mu \left( \Grad \vu + \Grad^t \vu - \frac{2}{3} \Div \vu \mathbb{I} \right) + \eta \Div \vu \mathbb{I}
$$
with constants $\mu >0$ and $\eta > 0$; and $s$ denotes the (specific) entropy.

The existence of finite energy weak solutions to (\ref{NSve}) for appropriately regular initial data has been established in \cite{MMMNPZ} under the assumption that the equations are coupled by the pressure in the form $$p(\vr,s) = \vr^\gamma \mathcal{T}(s), \ \  \ \gamma \geq \frac95,$$ where $\mathcal{T}(\cdot)$ is a given smooth and strictly monotone function from $R_+$ to $R_+$. In our case, we aim to take $\mathcal{T}(s) = \mbox{exp}((\gamma-1)s)$ with $\gamma= 1 + \frac1{c_v}$ to recover the state equations of a perfect gas
$$p(\vr, \vt) = \vr \vt \mbox{ and } s(\vr, \vt) = \log \left(\frac{\vt^{c_v}}\vr\right).$$  As we would like to have finite energy weak solutions to (\ref{NSve}) satisfying the {\it renormalized equation for entropy} for any $c_v >0$
we actually consider a pressure regularized "version" of (\ref{NSve}) instead and based on the results and know-how of \cite{MMMNPZ} we infer the existence of {\it entropy renormalized} weak solutions to   

\begin{equation} \label{asyst}
\begin{split} 
 \partial_t \vr + \Div (\vr \vu)&= 0, 
\\
\partial_t (\vr \vu) + \Div \left(\vr \vu\otimes \vu\right) + \Grad( \vr\vt) + \varepsilon \Grad (\vr\vt)^\frac\beta\gamma &= \ep \Div \mathbb{S}(\Grad\vu),
\\
\partial_t \int_\Omega \left[\frac 12 \vr |\bu|^2 + c_v \vr\vt+ \frac{\varepsilon}{\beta-1} (\vr\vt)^\frac\beta\gamma\right]\dx + \varepsilon \int_\Omega \mathbb{S}(\Grad \bu) : \Grad \bu  \, \dx
 &\le 0,
\\
 \partial_t \left(\vr \log \left(\frac{\vt^{c_v}}\vr\right)\right) + \Div \left(\vr \log \left(\frac{\vt^{c_v}}\vr\right)\vu\right)  &= 0,
\end{split}
\end{equation}
for some $\beta >>1$ (depending on $c_v>0$) and any $\ep >0$. We then take the limit $\ep \to 0$ based on  uniform a priori bounds on a sequence of solutions to \eqref{asyst} and  recover an rDMV solution to (\ref{p1}--\ref{pBC})  in the sense of Definition \ref{D1}. Actually, we get an rDMV solution that satisfies \eqref{p7} as an equality for any  continuous function $\chi$, see Theorem \ref{T1}.

\subsection{Existence of weak solutions to \eqref{asyst}} 

To show the existence of entropy renormalized solutions to \eqref{asyst} we 
follow the idea borrowed from \cite{MMMNPZ} and first consider  a pressure regularized isentropic Navier-Stokes system coupled with transport equation for entropy that describes the evolution of the mass density $\vr = \vr(t,x)$, the velocity field $\vu = \vu(t,x)$ and the pressure argument $Z=Z(t,x)$, in particular,

\begin{eqnarray}
\label{NS1} \partial_t \vr + \Div (\vr \vu)&=& 0, \\
\label{NS2} \partial_t (\vr \vu) + \Div \left(\vr \vu\otimes \vu\right) + \Grad Z^\gamma +  \delta\Grad Z^\beta &=&  \Div \mathbb{S}(\Grad\vu),
\\
\label{NS3} \partial_t Z + \Div \left(Z \vu \right)  &=& 0,
\end{eqnarray}
with the complete slip boundary conditions
\begin{equation} \label{uboundary_2}
\vu\cdot \vn|_{\partial \Omega} = 0 \mbox{ and } (\mathbb{S}(\Grad \vu)\cdot \vn) \times \vn|_{\partial \Omega}  =  0.
\end{equation}

 The system (\ref{NS1}--\ref{NS3}) without the pressure     "regularizing" term $\delta\Grad Z^\beta$ naturally appears in meteorology and astrophysics (see \cite {K}).  For  detailed discussion on existence of solutions see  \cite{MMMNPZ} and for singular limits see \cite{FKNZ}. In general, the idea is that $Z^\gamma = p$ is the pressure and under "enough" regularity we can 
relate it to the entropy $s$ and entropy transport equation
$$\partial_t (\vr s) + \Div \left(\vr s \vu \right)  = 0$$
 through  $p = \vr^\gamma \mathcal{T}(s)$.

The following result on existence of finite energy weak solutions to (\ref{NS1}--\ref{NS3}) can be found in \cite[Proposition 1]{MMMNPZ}.

\begin{Proposition} \label{prop1}
Let $\gamma >1$, $\beta \geq \max(\gamma,4)$ and $\delta>0$.
Then, given initial data $(\vr_0,Z_0,\bu_0) $ satisfying
\begin{equation} \label{reg_init}
\begin{split}
( \vr_0(\cdot),Z_0(\cdot), \vu_0(\cdot) ) \in C^{\infty}(\overline{\Omega},\R^5),\\
0 < c_\star \vr_0 \le Z_0 \le c^\star \vr_0 ~\text{in}~ \overline{\Omega}\mbox{ for some } 0< \ c_\star\leq c^\star <\infty,
\end{split}
\end{equation}
 there exists a finite energy weak solution $(\vr,Z,\bu)$ to problem (\ref{NS1})--(\ref{reg_init}) such that
\begin{equation*}
(\vr,Z,\bu) \in [ L^\infty(0,T;L^\beta(\Omega)) ]^2 \times L^2(0,T;W^{1,2}_0(\Omega,\R^3)),
\end{equation*}
\begin{equation}\label{inestep2delta}
0\le c_\star \vr \le Z \le c^\star \vr~\text{a.e in}~ (0,T)\times \Omega,
\end{equation}
and for any $\tau \in (0,T)$ we have:
\begin{description}
\item{(i)}
$ \vr \in C_w([0,T];L^\beta(\Omega))$ and the continuity equation (\ref{NS1}) is satisfied in the weak sense
\begin{equation}\label{contstep2delta}
\int_\Omega \vr(\tau,\cdot) \varphi(\tau,\cdot) \, \dx- \int_\Omega \vr_0 \varphi(0,\cdot) \, \dx= \int_0^\tau \int_\Omega \big(\vr \partial_t \varphi + \vr \bu \cdot \Grad \varphi\big) \, \dx\dt
\end{equation}
for any $\varphi \in C^1([0,T]\times \Ov{\Omega})$;
\item{(ii)}
 $ \vr \bu \in C_w([0,T];L^{\frac{2\beta}{\beta+1}}(\Omega,\R^3))$ and
the momentum equation (\ref{NS2}) is satisfied in the weak sense
\begin{multline}\label{momstep2delta}
\int_\Omega \vr \bu (\tau, \cdot) \cdot \bfphi(\tau,\cdot) \, \dx- \int_\Omega \vr_0 \vu_0 \cdot \bfphi(0,\cdot) \, \dx \\ = \int_0^\tau \int_\Omega \Big(\vr \vu \cdot \partial_ t \bfphi  + \vr \bu \otimes \bu : \Grad \bfphi + Z^\gamma \Div \bfphi   
+ \delta Z^\beta \Div \bfphi- \mathbb{S}(\Grad \bu) : \Grad \bfphi\Big) \, \dx\dt
\end{multline}
for any $ \bfphi \in C_c^1([0,T] \times \Ov{\Omega},\R^3)$, $\bfphi \cdot \vn|_{\partial \Omega} =0$;
\item{(iii)}
$ Z \in C_w([0,T];L^\beta(\Omega))$ and  equation (\ref{NS3}) is satisfied in the weak sense
\begin{equation}\label{entstep2delta}
\int_\Omega Z(\tau, \cdot) \varphi(\tau,\cdot) \, \dx- \int_\Omega Z_0 \varphi(0,\cdot) \, \dx= \int_0^\tau \int_\Omega \big(Z \partial_t \varphi + Z \bu \cdot \Grad \varphi\big) \, \dx\dt
\end{equation}
for any $\varphi \in C^1([0,T]\times\Ov{\Omega})$;
\item{(iv)}
the energy inequality
\begin{equation}\label{energystep2}
\begin{split}
\int_\Omega \left(\frac 12 \vr |\bu|^2 + \frac{1}{\gamma-1} Z^\gamma+ \frac{\delta}{\beta-1} Z^\beta\right)(\tau)\, \dx + \int_0^\tau \int_\Omega \mathbb{S}(\Grad \bu) : \Grad \bu  \, \dx\dt\\
 \le \int_\Omega \left(\frac 12 \vr_0 |\bu_0|^2 + \frac{1}{\gamma-1} Z_0^\gamma+ \frac{\delta}{\beta-1} Z_0^\beta\right)\, \dx 
\end{split}
\end{equation}
holds for a.a $ \tau \in (0,T) $.
\item{(v)}
 Moreover, equations  (\ref{NS1}), (\ref{NS3}) hold in the sense of renormalized solutions. That is, $(\vr,\bu, Z)$, extended by zero outside of $\Omega$, satisfy
\begin{equation}
\partial_t b(\vr) + \Div(b(\vr)\bu) + \big(b'(\vr)\vr-b(\vr)\big) \Div \bu = 0
\end{equation}
and 
\begin{equation}\label{renorent}
\partial_t b(Z) + \Div(b(Z)\bu) + \big(b'(Z)Z-b(Z)\big) \Div \bu = 0
\end{equation}
in  $\mathcal{D}'((0,T)\times \Omega)$ and $\mathcal{D}'((0,T)\times \R^3)$, where
$$
b \in C^1(\R), \quad b'(z)=0,\quad \forall z \in \R~\text{large enough.}
$$
\end{description}
\end{Proposition}

We note that \cite[Proposition 1]{MMMNPZ} gives the existence result under the no-slip boundary condition $\vu|_{(0,T)\times \partial \Omega} =0$, however, the result stays valid even for the complete slip conditions \eqref{uboundary_2} after obvious modification.

Now we show that a finite energy weak solution to (\ref{NS1}--\ref{uboundary_2}) obtained through Proposition \ref{prop1} is actually an entropy renormalized weak solution of \eqref{asyst}.

Let $c_v > 0$ and set $\gamma= 1 +\frac1{c_v}$ with $\beta \geq \mbox{max}\{\gamma,4\}$.  Then for any  $\delta >0$ and initial data $(\vr_0, Z_0,\vu_0)$ satisfying \eqref{reg_init} we get  from Proposition \ref{prop1} the existence of functions
$$(\vr,Z,\vu) \in  L^\infty(0,T;L^\beta(\Omega)) \times L^\infty(0,T; L^\beta(\Omega)) \times L^2(0,T;W^{1,2}_0(\Omega,\R^3))$$
that satisfy  (\ref{inestep2delta}--\ref{renorent}). Now thanks to \eqref{inestep2delta} and the regularity of $(\vr, Z)$ we
can define a.e. nonnegative function $\vt \in L^\infty(0,T;L^{\frac{\beta}{\gamma-1}}(\Omega))$ as
$$ \vt = \left\{\begin{array}{cl} \displaystyle\frac{Z^\gamma}\vr &\mbox{for } \vr > 0 \\ \\ 0& \mbox{for } \vr = 0 .\end{array} \right.$$
Moreover, following the same steps as in \cite[ Section 8.1]{MMMNPZ}  we can show (since $\beta \geq \frac95$) that 
$$c_* \leq \frac{Z}{\vr} \leq c^* \mbox{ a.e. in } (0,T) \times \Omega,$$
and  
$$\chi\left(\mathcal{T}^{-1}\left(\left(\frac{Z}\vr\right)^\gamma\right) \right)= \chi\left(  \frac1{\gamma-1} \log \left(\left(\frac{Z}\vr\right)^\gamma\right) \right) = \chi\left( \log \left( \frac{ \vt^{c_v} }{\vr} \right)\right)$$
satisfies the equation
\begin{equation} \label{NS8}
\begin{split}
\int_0^T \intO{ \left[ \vr \chi\left( \log \left( \frac{ \vt^{c_v} }{\vr} \right) \right)  \partial_t \varphi + \vr \chi \left( \log \left( \frac{ \vt^{c_v} }{\vr} \right) \right) \vu \cdot \Grad \varphi  \right] } \dt \\ = -
\intO{ \vr_0 \chi \left( \log \left( \frac{ \vt_0^{c_v} }{\vr_0} \right) \right) \varphi (0, \cdot) }
\end{split}
\end{equation}
for any $\varphi \in \DC([0,T) \times \Ov{\Omega})$, $\varphi \geq 0$ and any  $\chi \in C^1(R)$ and hence  by density argument for any $\chi \in C(R)$. Finally we note that $s =\log \left( \frac{ \vt^{c_v} }{\vr} \right) \in L^\infty((0,T)\times \Omega)$.

We have just showed the following existence result for \eqref{asyst}.

\begin{Proposition} \label{prop2}

Let $c_v >0$ and $\gamma=1 +\frac1{c_v}$. For any  $\beta \geq {\rm max}\{\gamma,4\}$,  $\ep >0$ and initial data $(\vr_0, \vt_0, \vu_0)$ satisfying 
\begin{equation} \label{asyst_init}
\begin{split}
( \vr_0(\cdot),\vt_0(\cdot),\vu_0(\cdot)) \in C^{\infty}(\overline{\Omega},\R^5), \\
 0 < c_\star \vr_0 \le \vt_0^{c_v}  \le c^\star \vr_0 ~\text{in}~ \overline{\Omega} \mbox{ for some } 0< \ c_\star\leq c^\star <\infty,
\end{split}
\end{equation}
there exists an entropy renormalized weak solution $(\vr,\vt,\vu)$ to \eqref{asyst} with \eqref{uboundary_2} and \eqref{asyst_init} such that
$$(\vr,\vt,\vu) \in  L^\infty(0,T;L^\beta(\Omega)) \times L^\infty(0,T; L^{\frac{\beta}{\gamma-1}}(\Omega)) \times L^2(0,T;W^{1,2}_0(\Omega,\R^3))$$ 
and  we have:
\begin{itemize}
\item
\begin{equation}\label{wf1}
 \int_0^T \int_\Omega \big(\vr \partial_t \varphi + \vr \bu \cdot \Grad \varphi\big) \, \dx\dt = 
- \int_\Omega \vr_0 \varphi(0,\cdot)
\end{equation}
for any $\varphi \in C_c^\infty([0,T)\times \Ov{\Omega})$;

\item
\begin{multline}\label{wf2}
\int_0^T \int_\Omega \Big(\vr \vu \cdot \partial_ t \bfphi  + \vr \bu \otimes \bu : \Grad \bfphi + \vr\vt \Div \bfphi   
+ \varepsilon (\vr\vt)^\frac\beta\gamma \Div \bfphi- \varepsilon\mathbb{S}(\Grad \bu) : \Grad \bfphi\Big) \, \dx\dt\\ =- \int_\Omega \vr_0 \vu_0 \cdot \bfphi(0,\cdot) \, \dx 
\end{multline}
for any $ \bfphi \in C_c^\infty([0,T) \times \Omega,\R^3)$, $\bfphi \cdot \vn |_{\partial \Omega} = 0$;

\item $\log \left( \frac{ \vt^{c_v} }{\vr} \right) \in L^\infty((0,T)\times \Omega)$ and 
\begin{equation} \label{wf3}
\begin{split}
\int_0^T \intO{ \left[ \vr \chi\left( \log \left( \frac{ \vt^{c_v} }{\vr} \right) \right)  \partial_t \varphi + \vr \chi \left( \log \left( \frac{ \vt^{c_v} }{\vr} \right) \right) \vu \cdot \Grad \varphi  \right] } \dt \\ = -
\intO{ \vr_0 \chi \left( \log \left( \frac{ \vt_0^{c_v} }{\vr_0} \right) \right) \varphi (0, \cdot) }
\end{split}
\end{equation}
for any $\varphi \in \DC([0,T) \times \Ov{\Omega})$, $\varphi \geq 0$ and any   $\chi \in C(R)$;

\item
\begin{equation}\label{wf4}
\begin{split}
\left[ \int_\Omega \left(\frac 12 \vr |\bu|^2 + c_v \vr\vt+ \frac{\varepsilon}{\beta-1} (\vr\vt)^\frac\beta\gamma \right)\, \dx\right]_{t=0}^{t=\tau} + \varepsilon\int_0^\tau \int_\Omega \mathbb{S}(\Grad \bu) : \Grad \bu  \, \dx\dt\leq 0
\end{split}
\end{equation}
holds for a.a $ \tau \in [0,T] $. 
\end{itemize}

\end{Proposition}

\subsection{Existence of rDMV solutions}

The goal of this section is to take a family of entropy renormalized weak solutions $(\vre, \vte, \vue)_{\ep > 0}$ to (\ref{wf1}--\ref{wf4})  with the initial data $(\vr_{0,\ep}, \vt_{0,\ep}, \vu_{0,\ep})_\ep$ generating a Young measure $U_{0}$ in an appropriate sense and show that as $\ep \to 0$ we recover an rDMV solution to the Euler system (\ref{p1}--\ref{pBC})  in the sense of Definition \ref{D1} with the  initial data $U_0$. A similar procedure has been done in \cite[Section 2.1]{BreFei17} or \cite[Section 3.5]{BreFei17A}.

To this end, we have to  discuss the following issues:
\begin{itemize}
\item Uniform bounds based on the energy estimates that will guarantee boundedness of the state variables $(\vre, \vme , E_\ep )$, where $\vme = \vre \vu_\ep$ and $E_\ep = \frac12 \vre|\vue|^2 + c_v \vre \vte$.
\item Showing that the viscosity term vanishes in the asymptotic limit.

\item Identifying the Young measure $\{ U_{t,x} \}$ associated to the family of (entropy renormalized) weak solutions.

\end{itemize}

Let the initial data $(\vr_{0,\ep}, \vt_{0,\ep}, \vu_{0,\ep})_\ep$ satisfy
\begin{equation} \label{initmeasgen1}
\begin{split}
( \vr_{0,\ep}(\cdot),\vt_{0,\ep}(\cdot),\vu_{0,\ep}(\cdot)) \in C^{\infty}(\overline{\Omega},\R^5),\\
0 < c_\star \vr_{0,\ep} \le \vt_{0,\ep}^{c_v} \le c^\star \vr_{0,\ep} ~\text{in}~ \overline{\Omega} \mbox{ for some } 0 < c_* \leq c^* < \infty,\\
 \intO{ \vr_{0, \ep} } \geq M_0 > 0, \\
\intO{ \left[ \frac{1}{2} \vr_{0, \ep}  |\vu_{0, \ep} |^2 + c_v \vr_{0, \ep} \vt_{0, \ep} + \frac{\ep}{\beta - 1}(\vr_{0,\ep} \vt_{0,\ep})^{\frac\beta\gamma}\right] } \leq e_0, \\
\end{split}
\end{equation}
 and generate a Young measure $U_{0}$ in the following sense
\begin{equation} \label{initmeasgen2}
\begin{split}
\intO{ \vr_{0,\ep} \phi } &\to \intO{ \left< U_{0,x}; \vr \right> \phi } \ \mbox{for any}\ \phi \in \DC(\Ov{\Omega}); \\
\intO{ \vr_{0,\ep} \vu_{0,\ep} \cdot \boldsymbol{\phi} } &\to \intO{ \left< U_{0,x}; \vm \right> \cdot \boldsymbol{\phi} } \\ 
&\mbox{for any}\ \boldsymbol{\phi} \in \DC(\Ov{\Omega}; R^3), \boldsymbol{\phi}  \cdot \vn |_{\partial \Omega} = 0;\\
\intO{ \left[ \frac{1}{2} \vr_{0,\ep} |\vu_{0,\ep}|^2 + c_v\vr_{0,\ep} \vt_{0,\ep} + \frac{\ep}{\beta - 1}(\vr_{0,\ep} \vt_{0,\ep})^{\frac\beta\gamma} \right] \phi }
&\to \intO{\left< U_{0,x}; E \right> \phi } \ \mbox{for any} \ \phi \in \DC(\Ov{\Omega});\\
\intO{  \vr_{0,\ep} \chi \left( \log \left( \frac{ \vt_{0,\ep}^{c_v}}{\vr_{0,\ep}}\right)\right)   \phi }
&\to \intO{ \left< U_{0,x}; \mathcal{S}_\chi(\vr, \vc{m}, E) \right> \phi } \\
 &\mbox{for any}\ \phi \in \DC(\Ov{\Omega}),\  \phi \geq 0,\mbox{ and any  }  \chi\in C(R).
\end{split}
\end{equation}

Then for any $\varepsilon >0$ there exists an entropy renormalized weak solution $(\vr_\ep, \vt_\ep, \vu_\ep)$  to (\ref{wf1}--\ref{wf4})  with the initial data $(\vr_{0,\ep}, \vt_{0,\ep}, \vu_{0,\ep})$ thanks to Proposition \ref{prop2}. For these solutions we deduce the following uniform bounds.

\subsubsection{Uniform bounds}

The total energy balance (\ref{wf4}) and \eqref{initmeasgen1} yield immediately
\begin{equation} \label{NS11}
{\rm ess} \sup_{t \in (0,T)} \intO{ \left[ E_\ep+  \frac{\varepsilon}{\beta-1} (\vre\vte)^\frac\beta\gamma\right]} + \varepsilon\int_0^T \int_\Omega \mathbb{S}(\Grad \vue) : \Grad \vue  \, \dx\dt \leq e_0.
\end{equation}
Further,  if we take an appropriate $\varphi = \psi(t)$, $\psi \geq 0$ in (\ref{wf3}) we get 
that 
\begin{equation*}
\intO{\vre(\tau, \cdot) \chi \left(\log \left(\frac{\vte^{c_v}(\tau, \cdot)}{\vre(\tau, \cdot)} \right)\right) }= \intO{\vr_{0,\ep}  \chi \left(\log \left(\frac{\vt_{0,\ep}^{c_v}}{\vr_{0,\ep}} \right)\right) } \mbox{ for a.a. } \tau \in (0,T),
\end{equation*}
which together with  the  continuous  function 
\[
 \chi(s)  \left\{ \begin{array}{lr} = 0 \ \mbox{for}\ \log c_* \leq s \leq \log c^*, \\ \\ <0 \ \mbox{ otherwise}, \end{array} \right.
\]
gives us the bounds 
\begin{equation} \label{NS12}
c_* \leq \frac{\vte^{c_v}(\tau, \cdot)}{\vre(\tau, \cdot)} \leq c_* \mbox{ whenever } 0< \vre
\end{equation}
and thus 
\begin{equation} \label{NS12a}
{\rm min}_{[\log  c_*, \log c^*]} \chi  \vre\leq \vre \chi\left( \log \left( \frac{ \vte^{c_v} }{\vre} \right) \right)  \leq{\rm max}_{[\log  c_*, \log c^*]} \chi  \vre \mbox{ a.e.  in } (0,T)\times \Omega 
\end{equation}
for any  $\chi \in C(R)$.

Moreover, using  (\ref{NS11}), \eqref{NS12} and $\vr \vu = \sqrt{\vr} \sqrt{\vr}\vu$ we conclude that
\begin{equation} \label{NS13}
{\rm ess} \sup_{t \in (0,T)} \| \vre(t, \cdot) \|_{L^{\gamma}(\Omega)} \leq c(e_0,c_*)
\end{equation}
and   
\begin{equation} \label{NS14}
{\rm ess} \sup_{t \in (0,T)} \| \vme (t, \cdot) \|_{L^{\frac{2\gamma}{\gamma + 1}}(\Omega)} \leq c(e_0,c_*),
\end{equation} 
where  $\frac{2\gamma}{\gamma + 1}> 1$.

Note that the estimates are uniform with respect to $\ep >0$ and strong enough to pass to the limit in the system (\ref{wf1}--\ref{wf4})  to generate an rDMV solution of the limit Euler system as soon as we show that the dissipative term in \eqref{wf2} vanishes in the asymptotic regime. However, that follows easily from \eqref{NS11}.

In view of the uniform bounds (\ref{NS11}), \eqref{NS13} and (\ref{NS14}) and the fundamental
theorem of the theory of Young measures (see e.g. Ball \cite{BALL2}), there is a subsequence of $(\vre, \vme, \Ee)_{\ep > 0}$ (not relabeled here) that generates
a Young measure $\{ U_{t,x} \}_{(t,x) \in (0,T) \times \Omega}$. Moreover, passing to the limit in the total energy balance (\ref{wf4}), we obtain
$$
\left[ \intO{ \left< U_{\tau,x} ; E \right> } \right]_{t = 0}^{t=\tau} + \mathcal{D} (\tau) = 0;
$$
for a.a. $\tau \in (0,T)$, where
\begin{equation} \label{ZD35}
\begin{split}
\mathcal{D}(\tau) &\geq \liminf_{\ep \to 0} \left(\intO{ \left[ \frac{1}{2} \vre |\vue|^2 +
c_v \vre \vte  + \frac{\ep}{\beta-1} (\vre \vte)^\frac\beta\gamma \right] (\tau, \cdot)} +\ep \int_0^\tau  \intO{ \mathbb{S}(\Grad \vu_\ep): \Grad \vu_\ep }\dt\right)\\ &- \intO{ \left< U_{\tau,x} ; E \right> } \ \mbox{for a.a}\ \tau \in (0,T).
\end{split}
\end{equation}

To pass to the limit in (\ref{wf1}) is straightforward and we obtain (\ref{p4}).

As already mentioned the dissipative term on the right-hand side  of the momentum balance  \eqref{wf2} vanishes and hence  we can pass to the limit to recover
(\ref{p5}). Note that the measure $\mu_C$ contains the concentration defect of the terms
\[
\vre \vue \otimes \vue,\quad  \vre\vte \ \ \mbox{ and }\ \  \ep (\vre \vte)^\frac\beta\gamma,
\]
and, by virtue of (\ref{ZD35}), it is controlled by $\mathcal{D}$ exactly as required in (\ref{p8}).

Finally, thanks to \eqref{NS12a} it is easy to perform the limit in the entropy balance
(\ref{wf3}) to obtain (\ref{p7}) even as an equality that holds for any  $\chi \in C(R)$.

We have shown the following result.

\begin{Theorem} \label{T1}

Let the transport coefficients satisfy
$\mu >0$ and  $\eta> 0$. Assume that  $(\vre, \vte, \vue)_{\ep > 0}$ is a family of finite energy weak solutions
to the Navier--Stokes system with entropy transport (\ref{asyst}) in the sense of Proposition \ref{prop2} and the initial data $(\vr_{0, \ep}, \vt_{0, \ep}, \vu_{0, \ep} )_{\ep > 0}$ satisfying \eqref{initmeasgen1} generate a Young measure $U_{0,x}$ in the sense of \eqref{initmeasgen2}.

Then (at least for a suitable subsequence)
\[
\left( \vre, \vre \vue, \frac{1}{2} \vre |\vue|^2 + c_v \vre \vte \right)_{\ep > 0}
\]
generates a Young measure $\{ U_{t,x} \}_{(t,x) \in (0,T) \times \Omega}$ that represents an rDMV solution of the Euler system (\ref{p1}), (\ref{p2})
in the sense of Definition \ref{D1}. Moreover, $ U_{t,x} $ satisfies the entropy balance
(\ref{p7}) as an equality that holds for any  $\chi \in C(R)$.
\end{Theorem}

\section{Generating rDMV solutions via Brenner's model}
\label{EB}

In this section we show that rDMV solutions can be obtained via the vanishing dissipation limit for
a two-velocity model introduced by H. Brenner in \cite{BREN2}--\cite{BREN1}. In particular, we consider the system

\begin{eqnarray}
\label{b1} \partial_t \vr + \Div (\vr \vu_m)&=& 0, \\
\label{b2} \partial_t (\vr \vu) + \Div \left(\vr \vu\otimes \vu_m\right) + \Grad p(\vr,\vt) &=&  \Div \mathbb{S}(\Grad \vu),
\end{eqnarray}
\begin{equation}
\begin{split}
\label{b3} \partial_t \left(\frac12\vr|\vu|^2 +  \vr e(\vr, \vt)\right) +
\Div \Big(\big(\frac12\vr|\vu|^2 +& \vr e (\vr, \vt) \big)\vu_m \Big)\\
+& \Div(p(\vr,\vt)\vu) + \Div \vq = \Div(\mathbb S(\Grad \vu)\cdot  \vu).
\end{split}
\end{equation}

The main idea of (\ref{b1}--\ref{b3}) is to introduce two velocity fields - $\vu$ and $\vu_m$ - interrelated through
\begin{equation}\label{b4}
\vu - \vu_m = K \Grad \log \vr,
\end{equation}
where $K \geq 0$ is a purely {\it phenomenological} coefficient. Here, it is  assumed that 
\begin{equation}\label{b5}
\vq = -\kappa \Grad \vt,
\end{equation}
\begin{equation}
\label{b6}
\S = \mu\left(\Grad \vu + \Grad^T \vu - \frac23 \Div \vu \I\right)
 + \eta \Div \vu \I
\end{equation}
are the Fourier heat flux and the viscous stress tensor, respectively.

As a matter of fact, Brenner's model has been thoroughly criticized and its relevance to fluid mechanics questioned in Oettinger et al. \cite{OeStLi}. On the other hand,
it is mathematically tractable and yields essentially better theory than the standard Navier--Stokes--Fourier system, see e.g. \cite{FeiVas}, Cai, Cao, Sun \cite{CaiCaoSun}.
Recently, the interest in ``two velocity models'' has been revived in Bresch et al. \cite{BrDeZa}, \cite{BrGiZa}.

Leaving apart the conceptual difficulties of the model, we show  that it generates in the vanishing dissipation limit an rDMV solution to the Euler system (\ref{p1}), (\ref{p2})
in the sense of Definition \ref{D1}. The crucial aspect of the analysis is a specific form of the coefficient $K$ in (\ref{b4}). Note that $K$ is taken constant in \cite{FeiVas} as well as in  Cai et al. \cite{CaiCaoSun}, while
Brenner proposed $K = \frac{\kappa}{c_p \vr}$, see \cite{BREN}, where $c_p$ denotes the specific heat at constant \emph{pressure}. On the other hand, Guermond and Popov \cite{GuePop} argue that the choice \begin{equation} \label{bconst}
K = \frac{\kappa}{c_v \vr},
\end{equation}
where  $c_v >0$ denotes the specific heat at constant \emph{volume}, should lead to an "ideal" numerical scheme  for approximating the complete Euler system.  During our analysis it turns out that \eqref{bconst} works well and hence we assume \eqref{bconst} together with the Boyle-Mariotte law 
\begin{equation}
\label{pg1}
p(\vr, \vt) = \vr \vt, \ \ \ e= c_v \vt.
\end{equation}

Under these circumstances, it is not difficult to show (see \cite{BreFei17A}) that the {\it renormalized entropy inequality} associated with (\ref{b1}--\ref{bconst}) reads  

\begin{equation} \label{b7}
\begin{split}
\partial_t &(\vr \chi(s) ) + \Div \left( \vr \chi(s) \vu_m \right) -
 \Div \left[ \frac{\kappa}{c_v} \Grad \chi(s) \right]  \\
&\geq \frac{\chi'(s)}{\vt} \mathbb{S}(\Grad \vu) : \Grad \vu +  {\kappa} \chi'(s) |\Grad \log{\vt}|^2  + \chi'(s) \frac{\kappa}{c_v} |\Grad \log \vr |^2
- \chi''(s) \frac{\kappa}{c_v} |\Grad s |^2
\end{split}
\end{equation}
for any non-decreasing function $\chi\in C^2(R)$, where $s= s(\vr, \vt)$ is the physical entropy of the system. In accordance with \eqref{pg1} we assume that the entropy $s$ is given by
\begin{equation}
\label{pg2}
 s(\vr, \vt) = c_v \log \vt - \log \vr.
\end{equation}

We further assume that the heat conductivity coefficient $\kappa$ and the viscosity coefficients $\mu$ and $\eta$ satisfy 
\begin{equation} \label{b8}
\kappa(\vr, \vt) \approx \vr, \ \ \ \mu(\vr, \vt) \approx \vr, \ \ \ \eta \equiv 0.
\end{equation}

The problem is closed by prescribing the complete slip and no-flux boundary conditions
\begin{equation} \label{bbc}
\vu\cdot \vn|_{\partial \Omega} = 0, \ \ \  (\mathbb{S}(\Grad \vu)\cdot \vn) \times \vn|_{\partial \Omega}  =  0, \ \ \ \Grad \vr \cdot \vn|_{\partial \Omega} = 0 \mbox{ and } \vq \cdot \vn |_{\partial \Omega} = 0.
\end{equation}

As the existence of weak or strong solutions to (\ref{b1}--\ref{bbc}) is not known at the moment we simply {\it assume} the existence of strong solutions here.

The main goal of this section is to show that a family  of classical solutions $(\vr_\ep, \vu_\ep, \vt_\ep)_{\ep >0}$ to 

\begin{equation}\label{bap4}
\vu - \vu_m = \ep \frac{\kappa}{c_v\vr} \Grad \log (\vr),
\end{equation}
\begin{eqnarray}
\label{bap1} \partial_t \vr + \Div (\vr \vu_m)&=& 0, \\
\label{bap2} \partial_t (\vr \vu) + \Div \left(\vr \vu\otimes \vu_m\right) + \Grad p(\vr,\vt) &=&  \ep\Div \mathbb{S}(\Grad \vu),
\end{eqnarray}
\begin{equation}
\begin{split}
\label{bap3} \partial_t \left(\frac12\vr|\vu|^2 +  \vr e(\vr, \vt)\right) +
\Div \Big(\big(\frac12\vr|\vu|^2 +& \vr e (\vr, \vt) \big)\vu_m \Big)\\
+& \Div(p(\vr,\vt)\vu) + \ep\Div \vq = \ep\Div(\mathbb S(\Grad \vu)\cdot  \vu),
\end{split}
\end{equation}
\begin{equation} \label{bap5}
\begin{split}
\partial_t &(\vr \chi(s) ) + \Div \left( \vr \chi(s) \vu_m \right) - \ep
 \Div \left[ \frac{\kappa}{c_v} \Grad \chi(s) \right]  \\
&= \ep\frac{\chi'(s)}{\vt} \mathbb{S}(\Grad \vu) : \Grad \vu +  \ep{\kappa} \chi'(s) |\Grad \log{\vt}|^2  + \ep\chi'(s) \frac{\kappa}{c_v} |\Grad \log \vr |^2
- \ep\chi''(s) \frac{\kappa}{c_v} |\Grad s |^2,
\end{split}
\end{equation}
for $\chi\in C^2(R)$, where $ s(\vr, \vt) = c_v \log \vt - \log \vr$,
generates an rDMV solution to the Euler system \eqref{p1} and \eqref{p2}. To this end, we first rewrite the system  (\ref{bap1}--\ref{bap5}) using \eqref{b5}, \eqref{pg1}, \eqref{bbc} and \eqref{bap4} as
\begin{eqnarray}
\label{nt1} \partial_t \vr + \Div (\vr \vu)&=& \frac{\ep}{c_v} \Div (\kappa \Grad \log \vr), \\
\label{nt2} \partial_t (\vr \vu) + \Div \left(\vr \vu\otimes \vu\right) - \frac{\ep}{c_v}\Div \left(\vu\otimes \kappa \Grad \log \vr\right) + \Grad (\vr\vt) &=&  \ep\Div \mathbb{S}(\Grad \vu),
\end{eqnarray}
\begin{equation}
\begin{split}
\label{nt3} \partial_t \intO{ \left(\frac12\vr|\vu|^2 +  c_v\vr \vt  \right)} = 0,
\end{split}
\end{equation}

\begin{equation} \label{nt4}
\begin{split}
\partial_t &(\vr \chi(s) ) + \Div \left( \vr \chi(s) \vu\right)  - \frac{\ep}{c_v}\Div \left( \vr \chi(s) \kappa \Grad \log \vr \right) - \frac{\ep}{c_v}
 \Div ( \kappa \Grad \chi(s) )  \\
&= \ep\frac{\chi'(s)}{\vt} \mathbb{S}(\Grad \vu) : \Grad \vu +  \ep{\kappa} \chi'(s) |\Grad \log{\vt}|^2  + \ep\chi'(s) \frac{\kappa}{c_v} |\Grad \log \vr |^2
- \ep\chi''(s) \frac{\kappa}{c_v} |\Grad s |^2,
\end{split}
\end{equation}
for $\chi\in C^2(R)$, where $ s(\vr, \vt) = c_v \log \vt - \log \vr$ and $\kappa$, $\mu$ and $\eta$ are given by \eqref{b8}. Second, we discuss the following issues:
\begin{itemize}
\item Uniform bounds based on the energy estimates that will guarantee boundedness of the state variables $(\vre, \vme , E_\ep )$, where $\vme = \vre \vu_\ep$ and $E_\ep = \frac12 \vre|\vue|^2 + c_v \vre \vte$.
\item Showing that the dissipation terms vanish in the asymptotic limit.

\item Identifying the Young measure $\{ U_{t,x} \}$ associated to the family of (renormalized) weak solutions.

\end{itemize}

Let the initial data $(\vr_{0,\ep}, \vt_{0,\ep}, \vu_{0,\ep})_\ep$ satisfy
\begin{equation} \label{initbren1}
\begin{split}
\vr_{0,\ep} >0, \ \ \   \intO{ \vr_{0, \ep} } \geq M_0 > 0, \ \ \ 
\vt_{0,\ep} >0, \ \ \ \log\left(\frac{\vt_{0,\ep}^{c_v}}{\vr_{0,\ep}}\right) \geq s_0 > -\infty, \\
\intO{ \left[ \frac{1}{2} \vr_{0, \ep}  |\vu_{0, \ep} |^2 + c_v \vr_{0, \ep} \vt_{0, \ep} \right] } \leq e_0, \\
\end{split}
\end{equation}
 and generate a Young measure $U_{0}$ in the following sense
\begin{equation} \label{initbren2}
\begin{split}
\intO{ \vr_{0,\ep} \phi } &\to \intO{ \left< U_{0,x}; \vr \right> \phi } \ \mbox{for any}\ \phi \in \DC(\Ov{\Omega}); \\
\intO{ \vr_{0,\ep} \vu_{0,\ep} \cdot \boldsymbol{\phi} } &\to \intO{ \left< U_{0,x}; \vm \right> \cdot \boldsymbol{\phi} } \\ 
&\mbox{for any}\ \boldsymbol{\phi} \in \DC(\Ov{\Omega}; R^3), \boldsymbol{\phi}  \cdot \vn |_{\partial \Omega} = 0;\\
\intO{ \left[ \frac{1}{2} \vr_{0,\ep} |\vu_{0,\ep}|^2 + c_v\vr_{0,\ep} \vt_{0,\ep}  \right] \phi }
&\to \intO{\left< U_{0,x}; E \right> \phi } \ \mbox{for any} \ \phi \in \DC(\Ov{\Omega});\\
\intO{  \vr_{0,\ep} \chi \left( \log \left( \frac{ \vt_{0,\ep}^{c_v}}{\vr_{0,\ep}}\right)\right)   \phi }
&\to \intO{ \left< U_{0,x}; \mathcal{S}_\chi(\vr, \vc{m}, E) \right> \phi } \\
 \mbox{for any}\ \phi \in \DC(\Ov{\Omega}),\  \phi \geq 0, \mbox{ and any }\chi &\mbox{ increasing concave satisfying } \chi(s) \leq \chi_\infty \mbox{ for all } s\in R.
\end{split}
\end{equation}

Assume that  for $\varepsilon >0$ the trio $(\vr_\ep, \vt_\ep, \vu_\ep)$ represents a classical solution  to (\ref{nt1}--\ref{nt4}) with \eqref{b8}, \eqref{bbc} and  the initial data $(\vr_{0,\ep}, \vt_{0,\ep}, \vu_{0,\ep})$ satisfying \eqref{initbren1} and \eqref{initbren2}. For these solutions we can deduce the following uniform bounds.

\subsection{Uniform bounds}

We get from the total energy balance \eqref{nt3} and \eqref{initbren1} the bound
\begin{equation} \label{bb1}
{\rm ess} \sup_{t \in (0,T)} \intO{ E_\ep(t,\cdot)} \leq e_0.
\end{equation}

For any $\chi$ such that $\chi' \geq 0$ (non-decreasing) and $\chi''\leq 0$ (concave) we get by integrating \eqref{nt4} over $\Omega$ and using \eqref{bbc} that 
\begin{equation} \label{ape2}
\partial_t \intO{\vre \chi(s_\ep)} \geq 0.
\end{equation}
In particular, taking $\chi$ in \eqref{ape2} such that 
\[
 \chi(s)  \left\{ \begin{array}{lr} \geq 0  \ \mbox{for}\ s_0 \leq s , \\ \\ <0 \ \mbox{for}\ s < s_0, \end{array} \right.
\]
 gives us together with \eqref{initbren1}  the {\it minimum principle}
\begin{equation*} 
s_0 \leq s_\ep(\tau, x) \mbox{ whenever } 0< \vre \mbox{ for a.a. } (\tau,x) \in (0,T)\times \Omega.
\end{equation*}
In other words, we have
$$ \vre \leq {\rm exp}(-s_0) \vte^{c_v}\mbox{ for a.a. } (\tau,x) \in (0,T)\times \Omega $$
and  hence from   (\ref{bb1}) it follows that
\begin{equation} \label{bb2}
{\rm ess} \sup_{t \in (0,T)} \| \vre(t, \cdot) \|_{L^{1+\frac1{c_v}}(\Omega)} \leq c(e_0,s_0).
\end{equation}

Moreover, using $\vr \vu = \sqrt{\vr} \sqrt{\vr}\vu$, \eqref{bb1} and \eqref{bb2} we conclude that
\begin{equation} \label{bb3}
{\rm ess} \sup_{t \in (0,T)} \| \vme (t, \cdot) \|_{L^{\frac{2c_v +2}{2c_v + 1}}(\Omega)} \leq c(e_0,s_0),
\end{equation} 
where  $\frac{2c_v +2}{2c_v + 1}>1$ for $c_v >0$. 

Finally, as
$$\vre |\log \vte|^p \leq \left\{\begin{array}{lr} c(s_0)\vte^{c_v} |\log \vte|^p \leq c(p,s_0) & \mbox{ if } \vte \leq 1 ,\\ \vre \vte & \mbox{ if } \vte \geq 1, \end{array}\right.$$
we get from \eqref{bb1} and \eqref{bb2} that 
\begin{equation}\label{bb4}
{\rm ess} \sup_{t \in (0,T)} \intO{ \vre(t,\cdot) \chi(s_\ep(t,\cdot))} \leq c(p,e_0,s_0)
\end{equation}
for any $\chi(s) \approx s^p$, $p \geq 1$ (independently of the $\chi$ class used elsewhere). Namely, we have it for $\chi(s) = s$ here.

On the other hand, if we  take  $\chi$ such that $\chi'(s) \geq c > 0$ for all $s \geq s_0$, $\chi''\leq 0$ (for example $\chi(s) = s$) we get (using \eqref{bb4}) by integrating \eqref{nt4} over $(0,T)\times\Omega$ that 
\begin{equation} \label{ape6}
\begin{split}
 \int_0^T\intO{\left(\ep\frac{1}{\vte} \mathbb{S}_\ep(\Grad \vue) : \Grad \vue +  \ep{\kappa_\ep}  |\Grad \log{\vte}|^2  + \ep \frac{1}{c_v} \kappa_\ep|\Grad \log \vre |^2 \right)
}\dt \leq c(e_0,s_0).
\end{split}
\end{equation}

\subsection{Estimates of dissipative terms}
We must show that the following terms 
$$\ep \kappa_\ep \Grad \log \vre, \ \ \ \ep\kappa_\ep(\vue \otimes \Grad \log \vre), \ \ \ \ep\mathbb{S}_\ep(\Grad \vue), \ \ \ \ep \kappa_\ep\chi(s_\ep) \Grad \log \vre, \ \ \ \ep \kappa_\ep \Grad \chi(s_\ep)$$
from \eqref{nt1}, \eqref{nt2} and \eqref{nt4} vanish as $\ep \to 0$. Since we have no control over them from energy inequality, we have to engage \eqref{ape6}. To that end we use the following estimates:
\begin{itemize}
\item
\begin{equation} \label{est1} \int_0^T \intO{ \ep \kappa_\ep \Grad \log \vre}\dt \leq \sqrt \ep (\|\kappa_\ep\|_{L^1((0,T)\times \Omega)} + \ep \|\sqrt{\kappa_\ep}\Grad \log \vre\|_{L^2((0,T)\times \Omega)}^2)
\end{equation}

\item
\begin{equation} \label{est2} \int_0^T \intO{ \ep\kappa_\ep(\vue \otimes \Grad \log \vre)}\dt \leq \sqrt{\ep}\big(  \|\sqrt{\kappa_\ep} \vue\|_{L^2((0,T)\times \Omega)}^2 + \ep \|\sqrt{\kappa_\ep}\Grad \log \vre\|_{L^2((0,T)\times \Omega)}^2\big)
\end{equation}

\item  
\begin{equation} \label{est3} \int_0^T \intO{\ep\mathbb{S}_\ep(\Grad \vue)}\dt \leq \sqrt{\ep} \Big( \|\vte \mu_\ep(\vre, \vte)\|_{L^1((0,T)\times \Omega)}
 + \int_0^T \intO{\ep \frac1\vte \mathbb{S}_\ep(\Grad \vue):\Grad \vue }\dt\Big)
\end{equation}

\item 
\begin{equation} \label{est4} \int_0^T \intO{ \ep \kappa_\ep \Grad \log \vte}\dt \leq \sqrt \ep (\|\kappa_\ep\|_{L^1((0,T)\times \Omega)} + \ep \|\sqrt{\kappa_\ep}\Grad \log \vte\|_{L^2((0,T)\times \Omega)}^2)
\end{equation}
\end{itemize}

\subsection{Vanishing dissipation limit}

In view of the uniform bounds \eqref{bb1}, \eqref{bb2}, \eqref{bb3} and the fundamental
theorem of the theory of Young measures (see e.g. Ball \cite{BALL2}), there is a subsequence of $(\vre, \vme, \Ee)_{\ep > 0}$ (not relabeled here) that generates
a Young measure $\{ U_{t,x} \}_{(t,x) \in (0,T) \times \Omega}$. Moreover, passing to the limit in the total energy balance (\ref{nt3}), we obtain
\begin{equation} \label{DL1}
\left[ \intO{ \left< U_{\tau,x} ; E \right> } \right]_{t = 0}^{t=\tau} + \mathcal{D} (\tau) = 0
\end{equation}
for a.a. $\tau \in (0,T)$, where
\begin{equation} \label{DL2}
\mathcal{D}(\tau) = \lim_{\ep \to 0} \intO{ \left[ \frac{1}{2} \vre |\vue|^2 +
c_v \vre\vte \right](\tau, \cdot) } - \intO{ \left< U_{\tau,x} ; E \right> } \ \mbox{for a.a}\ \tau \in (0,T).
\end{equation}

To pass to the limit in the weak formulation of (\ref{nt1}) to obtain (\ref{p4}) it is enough to observe that the term on the right-hand side tends to $0$ as $\ep \to 0$ thanks to \eqref{est1}, \eqref{ape6} and \eqref{bb2}.

We can pass to the limit in the weak formulation of the  momentum balance (\ref{nt2}) to recover
(\ref{p5}) since the extra terms vanish thanks to \eqref{est2}, \eqref{est3}, \eqref{ape6} and \eqref{bb1}. Note that the measure $\mu_C$ contains the concentration defect of the terms
\[
\vre \vue \otimes \vue\ \mbox{and}\ \vre \vte
\]
and, by virtue of (\ref{DL2}), it is controlled by $\mathcal{D}$ exactly as required in (\ref{p8}).

Finally, it remains to perform the limit in the entropy balance
(\ref{nt4}) to obtain (\ref{p7}). For $\chi$ such that $\chi'(s) \geq 0$ and $\chi''(s) \leq 0$ we see that  (\ref{nt4}) simplifies to

\begin{equation} \label{nt4s}
\partial_t (\vr \chi(s) ) + \Div \left( \vr \chi(s) \vu\right)  - \frac{\ep}{c_v}\Div \left( \vr \chi(s) \kappa \Grad \log \vr \right) - \frac{\ep}{c_v}
 \Div ( \kappa \Grad \chi(s) )  \geq 0.
\end{equation}
Hence for $\chi(s) \leq \chi_\infty $ for all $s \in R$ we can use the fact that $\chi'$ is bounded and  \eqref{est1}, \eqref{est4}, \eqref{ape6}, \eqref{bb2} to show that the dissipation terms vanish as $\ep \to 0$.

We have shown the following result:

\begin{Theorem} \label{T2}

Let the transport coefficients satisfy
$\mu \approx \vr$, $\eta \equiv 0$ and $\kappa \approx \vr$. Assume that  $(\vre, \vte, \vue)_{\ep > 0}$ is a family of classical solutions
to the Brenner's system (\ref{bap4}--\ref{bap5}) with the initial data $(\vr_{0, \ep}, \vt_{0, \ep}, \vu_{0, \ep} )_{\ep > 0}$ satisfying \eqref{initbren1} generate a Young measure $U_{0,x}$ in the sense of \eqref{initbren2}.

Then (at least for a suitable subsequence)
\[
\left( \vre, \vre \vue, \frac{1}{2} \vre |\vue|^2 + c_v \vre \vte \right)_{\ep > 0}
\]
generates a Young measure $\{ U_{t,x} \}_{(t,x) \in (0,T) \times \Omega}$ that represents an rDMV solution of the Euler system (\ref{p1}), (\ref{p2})
in the sense of Definition \ref{D1}. 
\end{Theorem}

\begin{Remark}
In our analysis we could also take $\mu (\vr, \vt) \approx   \frac1{\vt}$ or $\mu (\vr, \vt) \approx  \vr + \frac1{\vt}$. 
\end{Remark}

\section{Weak--strong uniqueness}
\label{WS}

The most notable property of rDMV solutions from Definition \ref{D1} is that they satisfy the  weak (measure-valued)--strong uniqueness principle, i.e., an rDMV solution and a strong solution starting from the same initial data coincide as long as the latter exists. In particular, in that case any sequence generating the  rDMV solution must converge to the strong solution and we can see rDMV solutions as a "tool" to prove it. This has already been applied to show the convergence of certain numerical schemes (see, e.g.,  \cite{FeiLuk}, \cite{FeiLukMiz}) and we hope to progress in that direction further.

To show the weak--strong uniqueness property for rDMV solutions we use the standard method of relative energy and relative energy inequality (see \cite{FeiNov10}). Namely, we use  that the difference between a strong solution and an rDMV solution emanating from the same initial data can be controlled by a coercive {\it relative energy} functional \eqref{WS2} and this functional satisfies a {\it relative energy inequality} \eqref{WS3}. Finally, using the fact that we work with a strong solution and an rDMV solution  we can show via the standard Gronwall argument that the relative energy functional is identically zero.

\subsection{Relative energy}

Let
$$
\tvr \in C^1([0,T] \times \Omega), \ \tvr > 0, \ \tvt \in C^1([0,T] \times \Omega), \ \tvt > 0, \ \tvu \in C^1([0,T] \times \Omega; R^3),
$$
be given. Following \cite{FeiNov10}, we introduce the \emph{ballistic free energy}
\[
H_{\tvt}(\vr, \vt) = c_v\vr \vt - \tvt \vr  \log\left(\frac{\vt^{c_v}}{\vr}\right),
\]
and the \emph{relative energy}
\[
\mathcal{E}\left(\vr, \vt, \vu \ \Big| \tvr, \tvt, \tvu \right) =
\frac{1}{2} \vr |\vu - \tvu|^2 + H_{\tvt}(\vr, \vt) -
\frac{ \partial H_{\tvt}(\tvr, \tvt) }{\partial \vr} (\vr - \tvr) - H_{\tvt} (\tvr, \tvt).
\]

We showed in \cite{BreFei17B} that $\mathcal{E}$  can be written in the conservative variables as
\begin{equation} \label{WS2}
\begin{split}
\mathcal{E} &\left( \vr, \vc{m}, E \Big| \tilde \vr,  \tilde{\vm}, \tilde E\right)
\\
&= - \tilde \vt \left[ \mathcal{S}(\vr,  \vc{m}, E) - \mathcal{S} (\tilde \vr,  \tilde{\vm}, \tilde E) \right. \\
&\ \ \ -\left. \partial_{\vr} \mathcal{S}(\tilde \vr,  \tilde{\vm}, \tilde E)(\vr - \tilde \vr)  - \nabla_{\vc{m}} \mathcal{S}(\tilde \vr, \tilde{\vm}, \tilde E) \cdot (\vc{m} - \tilde{\vc{m}}) -
\partial_E\mathcal{S}(\tilde \vr,  \tilde{\vm}, \tilde E )(E - \tilde E)
\right],
\end{split}
\end{equation}
where $\mathcal{S}$ is the {\it total entropy} without any normalization while  $\tilde {\vc{m}} = \tilde \vr \tilde \vu$ and   $\tilde{E} = \frac{1}{2} \tilde \vr |\tilde \vu|^2 + c_v \tilde \vr  \tilde \vt$.

It is interesting to see that the relative energy $\mathcal{E}$ is related to the relative entropy \` a la Dafermos \cite{Daf4} via a multiplicative factor proportional to the absolute temperature.

We showed in \cite{BreFei17} that rDMV solutions satisfy the following {\it relative energy inequality}
\begin{equation} \label{WS3}
\begin{split}
\Big[ \int_\Omega &\left< U_{t,x}; \mathcal{E} \left( \vr, \vc{m}, E \Big| \tilde \vr,  \tilde{\vm}, \tilde E \right) \right> \ \dx \Big]_{t = 0}^{t = \tau}
+ \intO{ \left[ \left< U_{0,x}; E \right> - \left< U_{\tau, x}; E \right> \right] }\\
\leq &- \int_0^\tau \intO{ \left[ \left< U_{t,x}; \mathcal{S}_\chi(\vr, \vc{m}, E)  \right> \partial_t \tilde \vt +
\left< U_{t,x};  \mathcal{S}_\chi (\vr, \vc{m}, E) \frac{\vc{m}}\vr \right> \cdot \Grad \tilde \vt \right] }\dt\\
&+ \int_0^\tau \intO{ \left[ \left< U_{t,x}; \vr \tvu - \vc{m} \right> \cdot \partial_t \tilde \vu  + \left< U_{t,x}; \frac{ (\vr \tvu - \vc{m}) \otimes \vc{m} }{\vr}
\right> : \Grad \tilde \vu  \right] } \dt \\ &- (\gamma - 1) \int_0^\tau \intO{ \left[ \left< U_{t,x}; E - \frac{1}{2} \frac{|\vc{m}|^2 }{\vr} \right> \Div \tilde \vu \right] }\dt \\
&+ \int_0^\tau \intO{ \left[ \left< U_{t,x}; \vr \right> \partial_t \tvt \log\left(\frac{\tvt^{c_v}}{\tvr}\right) + \left< U_{t,x}; \vc{m} \right> \cdot \Grad \tvt \log\left(\frac{\tvt^{c_v}}{\tvr}\right) \right] } \dt\\
&+ \int_0^\tau \intO{ \left[ \left< U_{t,x}; \tvr - \vr \right> \frac{1}{\tvr} \partial_t  (\tvr \tvt) - \left< U_{t,x}; \vc{m} \right> \cdot \frac{1}{\tvr} \Grad (\tvr \tvt) \right] } \dt\\
&+ \int_0^\tau \int_{\Omega} \Grad \tvu : {\rm d} \mu_C,
\end{split}
\end{equation}
with some suitably chosen fixed $\chi$  determined by $(\tvr, \tvt, \tvu)$ (see  \cite{BreFei17} for more details). We point out that the relation (\ref{WS3}) holds for \emph{any} trio of differentiable functions $(\tvr, \tvt, \tvu)$, $\tvr, \tvt > 0$.

\subsection{rDMV--strong uniqueness}

As a corollary of the relative energy inequality, we showed the weak--strong uniqueness principle in the class of rDMV solutions in \cite{BreFei17}.

The idea is to show that the terms on the right-hand side of (\ref{WS3}) can be absorbed by the time average of the left-hand side and hence by the means of standard Gronwall argument, the left-hand side must be identically zero on $(0,T)$. To this end, we use the coercivity properties of $\mathcal{E}$ following from (\ref{WS2}) and the fact that $\mathcal{S}$ is a concave function on its effective domain (see  \cite{BreFei17B}). We have obtained the following result in \cite{BreFei17}.

\begin{Theorem} \label{T3}

Let $c_v > 0$. Suppose that the Euler system (\ref{e1}--\ref{e3}) admits a continuously differentiable
solution $(\tvr, \tvt, \tvu)$ in $[0,T] \times \Omega$ emanating from the initial data
\[
\tvr_0 > 0,\ \tvt_0 > 0 \ \mbox{in}\ {\Omega}.
\]

Assume that  $\{ U_{t,x} \}_{(t,x) \in (0,T) \times \Omega}$ is an rDMV solution of the  system (\ref{p1}), (\ref{p2}) in the sense specified in Definition \ref{D1}, such that
\[
U_{0,x} = \delta_{\tvr_0(x), \tvr_0 \tvu_0 (x), \frac{1}{2} \tvr_0(x) |\tvu_0(x)|^2 + c_v \tvr_0 \tvt_0(x)} \ \mbox{for a.a.}\ x \in \Omega.
\]

Then
\[
U_{t,x} = \delta_{\tvr(t,x), \tvr \tvu (t,x), \frac{1}{2} \tvr(x) |\tvu(x)|^2 + c_v \tvr \tvt(t,x)} \ \mbox{for a.a.}\ (t,x) \in (0,T) \times \Omega.
\]
\end{Theorem}

Theorem \ref{T3} was given for periodic boundary conditions, that is on $\Omega = [0,1]^3|_{\{0,1\}}$, however it is easy to see that it stays valid even for the slip condition $\vu\cdot \vn|_{\partial \Omega} =0$.

In view of Theorems \ref{T1} (or \ref{T2}) and \ref{T3}, we immediately obtain the following corollary that can be seen as a version of the result in \cite{Fei2015A}.

\begin{Corollary} 

In addition to the hypotheses of Theorem \ref{T1} (or \ref{T2}) suppose that the limit Euler system (\ref{p1}), (\ref{p2}) admits a smooth ($C^1$) solution $(\vr, \vm, E)$ in
$[0,T] \times \Omega$.

Then
\[
\vre \to \vr,\ \vre \vue \to \vc{m}, \frac{1}{2} \vre |\vue|^2 + c_v \vre  \vte \to E \ \mbox{in}\ L^1((0,T) \times \Omega).
\]

\end{Corollary}

Indeed the fact that the limit DMV solution is represented by the Dirac masses implies (up to a subsequence) strong a.a. pointwise convergence. In addition,
the limit defect $\mathcal{D}$ vanishes which implies strong convergence in the $L^1-$norm.

\def\cprime{$'$} \def\ocirc#1{\ifmmode\setbox0=\hbox{$#1$}\dimen0=\ht0
  \advance\dimen0 by1pt\rlap{\hbox to\wd0{\hss\raise\dimen0
  \hbox{\hskip.2em$\scriptscriptstyle\circ$}\hss}}#1\else {\accent"17 #1}\fi}


\end{document}